\documentclass[twoside]{amsart}
\usepackage{amsmath}
\usepackage{amssymb}
\usepackage{enumerate}
\usepackage{amsthm}
\usepackage[colorlinks=true]{hyperref}

\theoremstyle{plain}
\newtheorem{theorem}{Theorem}[section]

\newtheorem{lemma}[theorem]{Lemma}
\newtheorem{corollary}[theorem]{Corollary}
\theoremstyle{definition}

\theoremstyle{remark}

\makeatletter
\newcommand{\maop}[1]{%
\ensuremath{\mathop{\operator@font #1}\nolimits}}
\makeatother
\newcommand{\downa}[1]{\ensuremath{\!\smash{\downarrow_{\!#1}}}}
\newcommand{\upa}[1]{\ensuremath{\!\smash{\uparrow^{\!#1}}}}

\newcommand{\Hom}{\mbox{\rm Hom}}

\newcommand{\End}{\mbox{\rm End}}

\newcommand{\Scal}{\mathcal{S}}

\newcommand{\NN}{\mathbb{N}}
\newcommand{\Acal}{\ensuremath{\mathcal{A}}}
\def\C/{\ensuremath{\mathbb C}}

\newcommand{\sgn}{\mathrm{sgn\,}}

\newcommand{\tScal}{\ensuremath{\tilde{\mathcal{S}}}}
\newcommand{\tAcal}{\ensuremath{\tilde{\mathcal{A}}}}
\newcommand{\tScalm}{\ensuremath{\tilde{\mathcal{S}}_{\{1,\ldots,m\}}}}
\newcommand{\tScalmm}{\ensuremath{\tilde{\mathcal{S}}_{\{m+1,\ldots,2m\}}}}
\newcommand{\Scalm}{\ensuremath{\mathcal{S}}_{\{1,\ldots,m\}}}
\newcommand{\Scalmm}{\ensuremath{\mathcal{S}}_{\{m+1,\ldots,2m\}}}
\newcommand{\ttau}{\ensuremath{\tilde{\tau}}}
\newcommand{\spinchar}[1]{\ensuremath{\smash <#1 \smash >}}
\newcommand{\Irr}{\mathrm{Irr}}
\newcommand{\Pcal}{\ensuremath{\mathcal{P}}}
\newcommand{\Dcal}{\ensuremath{\mathcal{D}}}
\newcommand{\Ocal}{\ensuremath{\mathcal{O}}}

\begin{document}

\title[Imprimitive Faithful Characters]{The Imprimitive Faithful 
Complex Characters of the Schur Covers of
the Symmetric and Alternating Groups}
\author{Daniel Nett}
\address{Lehrstuhl D f\"ur Mathematik, RWTH Aachen University,
52056 Aachen, Germany}
\email[D.N.]{Daniel.Nett@math.rwth-aachen.de}
\author{Felix Noeske}
\email[F.N.]{Felix.Noeske@math.rwth-aachen.de}

\thanks{Partially supported by the DFG grant HI 895/1-2}





\begin{abstract}
 Using combinatorics and character theory,
 we determine the imprimitive faithful complex characters, i.e., the
 irreducible faithful complex characters which are induced from proper
 subgroups, of the Schur covers of the symmetric and alternating groups. 
 Furthermore, for every imprimitive character we establish all 
 its minimal block stabilizers. As a corollary, we also determine 
 the monomial faithful characters of the Schur covers. 
\end{abstract}

\maketitle
\section{Introduction}\label{sec:intro}
In his seminal paper \cite{Aschbacher84}, \mbox{Aschbacher} gives a
subgroup structure theorem for the classical groups by defining eight
collections $C_1,\ldots,C_8$ of natural subgroups, as well as a class
$S$ of almost simple subgroups satisfying certain `irreducibility'
conditions. Building on that, in \cite{KleidmanLiebeck} Kleidman
and Liebeck have determined, given a finite almost simple classical
group $G$ of dimension at least 13, which members of $C_1,\ldots,C_8$
constitute maximal subgroups of $G$. This leaves only the maximal
members of $S$ to be determined, in order to obtain a classification of
the maximal subgroups of these classical groups.

One strategy is to determine those subgroups in $S$ which fail to be
maximal. By Aschbacher's theorem they are contained in maximal subgroups
which lie in one of the collections $C_2$, $C_4$, or $C_7$, or again in
$S$ (see \cite[\S 1.2]{KleidmanLiebeck}).

This note is a contribution to the analysis of the case $C_2$, as by
definition of this collection, its successful treatment amounts to
classifying all absolutely irreducible imprimitive representations
(respectively their characters) of the quasisimple groups and their
automorphism groups. Recall that an irreducible character $\chi$ of a
finite group $G$ is called \emph{imprimitive}, if it is induced from
a proper subgroup $H$ of $G$. In this case we call $H$ a \emph{block
stabilizer} of $\chi$.

In their articles \cite{DjokovicMalzanSn} and \cite{DjokovicMalzanAn},
Djokovi\'c and Malzan have already determined all imprimitive ordinary
characters of the symmetric and the alternating groups, allowing us to
focus on the faithful characters of the covers of these groups.
The Schur covers of the symmetric groups are non-split central extensions
\begin{equation*}
    1\longrightarrow Z\longrightarrow \tScal_n 
    \overset{\theta}{\longrightarrow} S_n \longrightarrow 1,
\end{equation*}
for which the Schur multiplier $Z$ has order two, and $Z \le \tScal_n'$,
provided $n$ is at least four. For $n<4$ the Schur multiplier is
trivial, so let $n\ge 4$ hold throughout the present paper. There are
two isoclinic Schur covers of $\Scal_n$ which are only isomorphic if
$n=6$. In this note we adhere to Schur's choice (see \cite{Schur1911})
and consider the group defined by the presentation
\begin{align*}
\begin{split}\label{praestilde}
\tScal_n \cong \left< z,t_1,\ldots,t_{n-1}|\right. &
z^2=1, t_i^2=z, \: 1\le i \le n-1,\\
& (t_it_{i+1})^3=z, \: 1\le i \le n-2,\\
&\left. t_it_j=zt_jt_i, \: 1\le i < j \le n-1,|i-j|>1 \right>, 
\end{split}
\end{align*}
i.e., transpositions lift to order four. Therefore, given the presentation
\begin{align*}
\begin{split}
\Scal_n \cong \left< s_1,\ldots,s_{n-1}|\right. &
 s_i^2=1, \: 1\le i \le n-1,\\
& (s_is_{i+1})^3=1, \: 1\le i \le n-2,\\
&\left. s_is_j=s_js_i, \: 1\le i < j \le n-1,|i-j|>1 \right> 
\end{split}
\end{align*}
of $\Scal_n$, the epimorphism $\theta$ is defined by mapping $s_i$ to
$t_i$ for all $1\le i \le n-1$. For any $g \in \Scal_n$ we define its
\emph{standard lift} to $\tScal_n$ to be the element which is obtained
by replacing $s_i$ in the product expansion of $g$ into transpositions
by $t_i$ for every $1\le i \le n-1$.

Note that even though we are focusing on $\tScal_n$, we also
obtain information on the imprimitive characters of the isoclinic
group $\hat{\Scal}_n$: By \cite[Section 6.7]{ATLAS} (see also
\cite[III.5]{BeylTappe} for details) the characters of $\tScal_n$ may be
tensored with an appropriate linear character to give the characters of
the isoclinic cover $\hat{\Scal}_n$. Hence if a character of $\tScal_n$
is imprimitive, so is its corresponding character of $\hat{\Scal}_n$.

Except for $n=6$ or $n=7$ the Schur covers of the alternating groups
are simply given by the derived subgroups of the Schur covers of the
symmetric groups. As isoclinic groups possess isomorphic derived
subgroups, there is always only one cover up to isomorphism. Hence we
may define $\tAcal_n := \tScal_n'$.
When $n=6$ or $n=7$, the Schur multiplier of $\Acal_n$ is exceptional:
In these cases it is cyclic of order six, giving rise to the two
non-split extensions $6\dot{}\Acal_6$ and $6\dot{}\Acal_7$.

Before we give the main result of this paper, let us briefly review
some basic facts of the faithful complex characters of $\tScal_n$ to
introduce some notation in the process: Of course, when dealing with
symmetric groups, the proper combinatorial objects to consider are
partitions of $n$, whose set we denote by $\Pcal_n$, and whose parts we
assume to be ordered ascendingly.
Furthermore, if $\pi \in \Pcal_n$ is the cycle type of some $g\in
\Scal_n$, then we set $\sigma(g):= n - \ell(\pi) \pmod 2$, where
$\ell(\pi)$ is the length of $\pi$. Hence, $(-1)^{\sigma(g)}$ gives the
familiar sign character $\sgn$, and we extend $\sigma$ to $\tScal_n$
by setting $\sigma(\tilde{g}) := \sigma( \theta(\tilde{g}) )$ for
$\tilde{g} \in \tScal_n$.
Also note that, as $\tAcal_n$ is the full $\theta$-preimage of the
alternating group $\Acal_n$, we have $\ker \sigma = \tAcal_n$. 
For any
subgroup $G$ of $\tScal_n$ containing the center $Z$
the set of its ordinary characters faithful on $Z$
is denoted by $\Irr_-(G)$, and the elements of $\Irr_-(\tScal_n)$
and $\Irr_-(\tAcal_n)$ are called \emph{spin characters}. 
The spin
characters of $\tScal_n$ are parameterized by strict partitions, i.e., by
partitions of $n$ all of whose parts are distinct. We write $\Dcal_n$
for the set of strict partitions of $n$, $\Dcal_n^+$ for all even
strict partitions, and $\Dcal_n^-$ for all strict partitions which
are odd. Given $\lambda \in \Dcal_n$ we denote the corresponding spin
character by $\spinchar{\lambda}$. If $\lambda$ is an odd partition, it
gives rise to two \emph{associate} characters $\spinchar{\lambda}$ and
$\spinchar{\lambda}^a:= \sgn \otimes \spinchar{\lambda}$.

With a minimum of notation in place, we can now state our result on the
imprimitive faithful ordinary characters of $\tScal_n$ and $\tAcal_n$.
We do so by giving all triples $(H,\varphi,\chi)$, where $\chi$ is an
imprimitive faithful character of $G \in \{\tScal_n,\tAcal_n\}$ such
that $\chi=\varphi\upa{G}$ for some $\varphi\in\Irr_-(H)$ and $H\leq G$
is a subgroup minimal with this property, i.e., $H$ does not properly
contain a block stabilizer of $\chi$.
In the following, we refer to such triples as \emph{minimal triples}.
\begin{theorem} \label{thm:Sn} \label{thm:mainSn} \label{la:verticesSn} 
  For $\tScal_n$ the minimal triples $(H,\varphi,\chi)$ 
  are the following:
  \begin{itemize}
    \item[(i)] $H=\tAcal_n$, $\varphi$ is a constituent of $\spinchar{\lambda}
      \downa{\tAcal_n}$ for some $\lambda\in\Dcal_n^+$ and $\chi=
      \spinchar\lambda$, except when $n=6$ and $\lambda=(4,2)$, or $n=9$ and
      $\lambda=(6,2,1)$.
    \item[(ii)] $H=\tScal_{n-1}$, $\varphi=\spinchar\mu$ for 
     $\mu=(m,m-1,\dots,1)\in \Dcal_{n-1}^-$ 
     with $m \equiv 2,3 \pmod 4$, and 
      $\chi=\spinchar\lambda$ with
      $\lambda=(m+1,m-1,m-2,\dots,1)\in\Dcal_n^+$.
    \item[(iii)] For $n=6$ we have $H=3^2\colon 8$, $\varphi$ is an extension of
     either 
     linear character of order four of the subgroup $3^2\colon 4$,
     and $\chi = \spinchar{4,2}$ is the unique
     spin character of degree 20. 
    \item[(iv)] For $n=9$: $H=2\times L_2(8)\colon3$, 
     $\varphi$ is 
     a linear character of order six, and $\chi=\spinchar{6,2,1}$ is the 
     unique spin character of degree 240.
  \end{itemize}
\end{theorem}

\begin{theorem} \label{thm:An}
 For $\tAcal_n$ the minimal triples $(H,\varphi,\chi)$ are as follows:
  \begin{itemize}
    \item[(i)] $H=\tAcal_{n-1}$, $\varphi$ is either constituent of 
     $\spinchar\mu
      \downa{\tAcal_{n-1}}$ for $\mu=(m,m-1,\dots,1)\in\Dcal_{n-1}^+$ and
      $\chi=\spinchar{\lambda}\downa{\tAcal_n}$ with $\lambda = 
      (m+1,m-1,\dots,1)\in \Dcal_n^-$.
    \item[(ii)] $n=6$: $H=3^2\colon 8$, $\varphi$ is an extension of either
      linear character of order four of the subgroup $3^2\colon4$.
    Both linear characters have a pair of extensions
    in $H$. The members of each pair induce to the constituents of
    $\spinchar{4,2}\downa{\tAcal_6}$, i.e., the two faithful characters of
    degree 10.
    \item[(iii)] $n=9$: $H=2\times L_2(8)\colon3$, 
     taking $\varphi$ to be 
     one of the irreducible linear characters of order six of $H$ 
     gives either of the two faithful
     characters of degree 120, i.e., $\chi$ is a constituent of
     $\spinchar{6,2,1}\downa{\tAcal_9}$.
  \end{itemize}
\end{theorem}

The proofs of Theorems~\ref{thm:Sn} and \ref{thm:An} are given in
Sections~\ref{sec:symm} and \ref{sec:alt}. At the end of
Section~\ref{sec:alt} we also classify the imprimitive faithful characters
of the exceptions $6\dot{}\Acal_6$, $6\dot{}\Acal_7$, $3\dot{}\Acal_6$ and
$3\dot{}\Acal_7$ (see Theorem~\ref{thm:6An}).

Our strategy is as follows: We begin in Section~\ref{sec:reductions}
by determining the proper subgroups of $\tScal_n$ and $\tAcal_n$
which are viable candidates for block stabilizers. Of course, by
transitivity of induction, it is sufficient to restrict our attention
to maximal subgroups. As these are $\theta$-preimages of maximal
subgroups of $\tScal_n$ or $\tAcal_n$, they form three classes of
subgroups, distinguished by the type of the natural action of their
image under $\theta$: We analyze intransitive, imprimitive and
primitive subgroups separately. For the first two classes this analysis
is straightforward, shortening our list of viable candidates
considerably for the imprimitive subgroups. For the primitive subgroups,
whose classification in general is still open, we make use of a
classification by \mbox{Kleidman} and Wales in \cite{KleidmanWales91},
who determined the primitive subgroups of $\Scal_n$ whose order is at
least $2^{n-4}$. 

In Section~\ref{sec:comb} we lay the combinatorial foundations
needed to successfully tackle the double covers of the symmetric
groups: The main tools are the Branching Rule and an analogue of the
Littlewood-Richardson rule for spin characters due to Stembridge (see
\cite{Stembridge89}). The latter provides a means to describe the
constituents of a faithful character induced from a maximal intransitive
subgroup of $\tScal_n$ (Stembridge refers to such a character as a 
\emph{projective outer product}).

With the combinatorics in place, in Section~\ref{sec:symm} we complete
Bessenrodt's classification of the projective outer products which
are multiplicity free (see \cite{Bessenrodt02}) with a result of the
first author's diploma thesis \cite{Nett07}. As a consequence we obtain
the imprimitive characters of $\tScal_n$ whose block stabilizers are
maximal intransitive subgroups. Furthermore, this information is used to
analyze the situation for imprimitive subgroups of $\tScal_n$: Together
with character theoretic arguments and Clifford Theory we are able to
classify the imprimitive characters induced from maximal imprimitive
subgroups.

In Section~\ref{sec:alt} we employ Clifford Theory again to apply
our results of the previous section to the double covers of the
alternating groups. Lastly, we deal with the exceptional Schur covers
$6\dot{}\Acal_6$ and $6\dot{}\Acal_7$ with the help of the ATLAS
\cite{ATLAS} and some \textsf{GAP}-calculations (see \cite{GAP4}). The
triple covers $3\dot{}\Acal_6$ and $3\dot{}\Acal_7$ are also considered.

\section{Reductions}\label{sec:reductions}
In order to determine the imprimitive irreducible characters of
any finite group $G$, it is sufficient to consider the irreducible
characters of its maximal subgroups, and single out those which induce
irreducibly. To this end let us begin by stating two results which allow
us to narrow down the list of subgroup candidates for block stabilizers
in $G$. First, by Lemma~\ref{lem:sqrt} the orders of the candidates must
not be too small.
\begin{lemma}
 Let $G$ be a finite group, $H \le G$ and $\chi \in \Irr(H)$. If
 $|G:H|^2 > |G|$, then $\chi\upa{G}$ is reducible. \label{lem:sqrt}
\end{lemma}
\begin{proof}
 Since $\chi\upa{G}(1) = \chi(1)|G:H|$, the hypothesis forces
 $\chi\upa{G}(1)^2 > |G|$, hence $\chi\upa{G}$ cannot be irreducible.
\end{proof}
Secondly, the following consequence of Mackey's Theorem, which is
already used implicitly in \cite{DjokovicMalzanSn, DjokovicMalzanAn},
enables us to eliminate candidates, too.
\begin{lemma}\label{lem:centralizing_intersection}
 Let $G$ be a finite group, $H \le G$ and $\chi \in \Irr(H)$. If there
 exists a non-trivial element $t \in G$ which centralizes $H^t \cap H$,
 then $\chi\upa{G}$ is reducible.
\end{lemma}
\begin{proof}
 Let $M$ be the $\C/H$-module affording $\chi$, and consider the
 endomorphism ring $E:=\End_{\C/G}(M\upa{G})$. By adjointness and
 Mackey's Theorem we have \[ \dim E = \dim \left( \bigoplus_{y
 \in H\backslash G /H} \Hom_{\C/(H^y\cap H)}( M\downa{H^y\cap H},
 M^y\downa{H^y \cap H}) \right).\] By the hypothesis $M^t\downa{H^t\cap
 H} \cong M\downa{H^t\cap H}$, thus taking one of the double coset
 representatives to be $t$, we obtain
 \begin{eqnarray*}
  \dim E &=& \sum_{y \in H \backslash G /H} \dim \Hom_{\C/(H^y\cap H)}(M
  \downa{H^y \cap H}, M^y\downa{H^y \cap H})\\ &\ge& \dim \End_{\C/H}(M)
  + \dim \End_{\C/(H^t\cap H)}(M\downa{H^t \cap H})\\ &\ge& 2.
 \end{eqnarray*}
\end{proof}

From now on let $G = \tScal_n$ or $G= \tAcal_n$. As the maximal
subgroups of $G$ are the preimages under $\theta$ of the maximal
subgroups of $\theta(G)$, we construct our list of viable candidates by
considering subgroups of the symmetric and alternating groups first.

The subgroups of $\Scal_n$ and $\Acal_n$ fall into three main classes
distinguishable by the type of their natural action on the set
$\{1,\ldots,n\}$: The intransitive subgroups, the subgroups which are
imprimitive and transitive, and those which act primitively.

The intransitive maximal subgroups of $\Scal_n$ are the maximal
parabolic subgroups isomorphic to a group of the form $\Scal_l\times
\Scal_{n-l}$ for $1\le l \le \lfloor \frac{n}{2}\rfloor$. As every
intransitive subgroup of $\Acal_n$ is contained in an intransitive
maximal subgroup of $\Scal_n$, the intransitive maximal subgroups of
$\Acal_n$ are given by $(\Scal_l\times \Scal_{n-l}) \cap \Acal_n$. Hence
we obtain the following lemma.
\begin{lemma}\label{lem:intransitive}
 Let $H$ be a maximal subgroup of $G$ such that $\theta(H)$ acts
 intransitively, then $H$ is the preimage of $(\Scal_l\times
 \Scal_{n-l})\cap \theta(G)$.
\end{lemma}

The imprimitive and transitive maximal subgroups of $\Scal_n$ are of
the form $\Scal_{n/k} \wr \Scal_k$ for $1\le k\le \lfloor \frac{n}{2}
\rfloor$ with $k \:|\: n$. Again as every imprimitive and transitive
subgroup of $\Acal_n$ is contained in an imprimitive and transitive
maximal subgroup of $\Scal_n$, the imprimitive and transitive maximal
subgroups of $\Acal_n$ have the form $(\Scal_{n/k} \wr \Scal_k)\cap
\Acal_n$, where $k$ is as above. However, not all of these maximal
subgroups are viable candidates, as we can rule out the majority of
cases with the help of Lemma~\ref{lem:centralizing_intersection}.
\begin{lemma}
 If $H$ is the block stabilizer of some irreducible and imprimitive
 character $\chi \in \Irr(G)$, and $H$ is the preimage of a transitive
 and imprimitive maximal subgroup of $G$, then $n$ is even and $H$
 is the preimage of $(\Scal_{n/2}\wr \Scal_2) \cap \theta(G)$.
 \label{lem:transitive_imprimitive}
\end{lemma}
\begin{proof}
 Assume $n=kl$ and let $H$ be the preimage of $(\Scal_{l}\wr
 \Scal_k)\cap \theta(G)$ for some $k\ge 3$. Without loss of generality
 let $B_i:=\{(i-1)l+1,\ldots,il\}$ for $i=1,\ldots,k$ be the blocks
 of the action of $\theta(H)$ on $\{1,\ldots,n\}$. Furthermore,
 let $t:=(1, l+1, 2l+1) \in \theta(G)$. 
 Then the centralizer of
 $\theta(H)^t \cap \theta(H)$ in $\theta(G)$ 
 contains $t$. As the order
 of $t$ is odd, there exists an element $\tilde{t} \in \theta^{-1}(t)$
 which is of odd order, too. Now, for any $h \in C_{G}(H^{\tilde{t}}\cap
 H)$ we have $\theta(h^{\tilde{t}}) = \theta(h)^t = \theta(h)$, and
 hence $h^{\tilde{t}} \in \{ h, zh\}$. By our choice of $\tilde{t}$,
 we conclude $h^{\tilde{t}} = h$, i.e., the element $\tilde{t}$
 centralizes $C_{G}(H^{\tilde{t}}\cap H)$. Therefore
 Lemma~\ref{lem:centralizing_intersection} gives the claim.
\end{proof}

In the case of a primitive action,
for our
treatment when $n\ge 5$ we employ a classification of the primitive
subgroups of $\Scal_n$ which do not contain $\Acal_n$, and whose order 
is at least $2^{n-4}$. This is provided by Kleidman and Wales in
\cite[Proposition~6.2]{KleidmanWales91} on the basis of the O'Nan-Scott
Theorem (see, for example, \cite{LiebeckPraegerSaxl}).
Together with Lemma~\ref{lem:sqrt} we can derive the
following result on irreducible characters 
induced from primitive subgroups.
\begin{lemma}\label{la:kleidmanwales}
 Let $H\neq \tAcal_n$ be the preimage of a primitive subgroup of
 $G$ such that $H$ is the block stabilizer of some irreducible
 and imprimitive $\chi \in \Irr_-(G)$. Then $n=9$ and $H\cong
 2\times\mathrm{L}_2(8)\colon3$. Let $\phi$ and $\phi'$ denote the two
 linear irreducible characters of order six of $H$. For $G=\tAcal_9$
 we obtain $\chi = \phi\upa{G}$ or $\chi = \phi'\upa{G}$ and for
 $G=\tScal_9$ we obtain $\chi = \phi\upa{G} = \phi'\upa{G}$.
 \label{lem:primitive}
\end{lemma}
\begin{proof}
 If $n=4$ then it is easily verified that no maximal subgroup of $G$
 fulfilling the hypothesis is primitive. Therefore let $n \ge 5$,
 and we commence by considering the case $|\theta(H)| \ge 2^{n-4}$.
 Let $G = \tScal_n$. With the help of Lemma~\ref{lem:sqrt} we can
 work through the list of \cite[Proposition 6.2]{KleidmanWales91}
 and verify for most groups that they are not block stabilizers of
 some $\chi \in \Irr_-(G)$. There remains a small list of groups for
 which we have to take a closer look. These are the preimages of
 $\Acal_5$ and $\Scal_5$ for $n=6$, $\mathrm{L}_2(7)$ for $n=7$,
 $2^3\colon\mathrm{L}_3(2)$ and $\mathrm{L}_3(2).2$ for $n=8$,
 $3^2\colon\tScal_4$ and $\mathrm{L}_2(8)\colon 3$ for $n=9$,
 $\mathsf{M}_{11}$ for $n=11$, and finally $\mathsf{M}_{12}$ for
 $n=12$. Except in the case $\mathrm{L}_2(8)\colon 3$ for $n=9$, the
 character table of $G$ shows that there is no spin character whose
 degree is divisible by the index of the given groups in $\Scal_n$
 (which of course is the index of the respective preimage in $G$).
 The preimage of $\mathrm{L}_2(8)\colon 3$ is the subgroup $2\times
 \mathrm{L}_2(8)\colon 3 \le G$. Its commutator factor group is cyclic
 of order six, and with the aid of $\textsf{GAP}$ we verify that the
 inflation of its two faithful characters induce irreducibly to give the
 same character of $G$.

 Let $G=\tAcal_n$. As the primitive subgroups of $\Acal_n$ are obtained
 by intersection those of $\Scal_n$ with $\Acal_n$, we may argue as in
 the case $G=\tScal_n$ by going through the corresponding subgroups of
 $\Acal_n$ using \cite[Proposition 6.2]{KleidmanWales91}. Again only
 the subgroup $2\times \mathrm{L}_2(8)\colon 3$ which is maximal in
 $\tAcal_9$ gives rise to two irreducible and imprimitive characters
 of $\tAcal_9$: Here both inflated faithful characters of the
 commutator factor group induce to two non-isomorphic characters.

 On the other hand, if $|H| < 2^{n-4}$, then as $|H|^2 < |G|$ for both
 $G=\tScal_n$ and $G= \tAcal_n$, we have $|G:H|^2 > |G|$, and no such
 group is a block stabilizer of some irreducible and imprimitive $\chi
 \in \Irr(G)$ by Lemma~\ref{lem:sqrt}.
\end{proof}

If $G=\tScal_n$ Lemma~\ref{la:kleidmanwales} leaves the case
$H=\tAcal_n$, i.e., it remains to decide which characters of
$\tScal_n$ are induced from $\tAcal_n$. But this is well known (see,
for example, \cite[Theorem~4.2]{HoffmanHumphreys}): A spin character
$\spinchar{\lambda} \in \Irr_-(\tScal_n)$ is induced from $\tAcal_n$
if and only if $\lambda$ is an even partition. Therefore we will not
concern ourselves with it any further in the sequel.

\section{Combinatorics}\label{sec:comb}
In this section we establish the combinatorial framework needed in
Sections~\ref{sec:symm} and \ref{sec:alt}. We fix the notation and for
the convenience of the reader we collect the necessary theorems on which
further proofs are based.

One of them is the Branching Rule for spin characters. To state it and
to provide further notation which will be useful for some arguments
used in Sections~\ref{sec:symm} and \ref{sec:alt}, we have to describe
possible enlargements of partitions. Therefore for a strict
partition $\lambda$
of $n$, let $N(\lambda)$ be the subset of $\Dcal_{n+1}$ consisting of
all partitions which can be obtained by adding 1 to one of the parts of
$\lambda$. If $\lambda$ does not contain 1 as a part, let $\lambda^+:=
\lambda \cup (1)$, where for two partitions $\lambda$ and $\mu$ we
denote by $\lambda \cup \mu$ the partition obtained by forming the union
of both sets of parts. Note that $\lambda^+$ is not in $N(\lambda)$.
Furthermore defining
\begin{equation*}
 \spinchar{\mu}^{*}=
 \begin{cases}
  \spinchar{\mu}, &\text{if }\mu\text{ is even},\\
  \spinchar{\mu} + \spinchar{\mu}^{a}, &\text{if }\mu\text{ is odd}.
 \end{cases}
\end{equation*}
for a strict partition $\mu$, allows us to state the Branching Rule for
inducing spin characters from $\tScal_n$ to $\tScal_{n+1}$.
\begin{theorem}[{\cite[Theorem 10.2]{HoffmanHumphreys}}]
\label{thm:BR}
 Let $\lambda$ be a strict partition of $n$. Then
 \[
  \spinchar{\lambda}\upa{\tilde S_{n+1}} = 
  (1- \delta_{1,\lambda_{\ell(\lambda)}}) \spinchar{\lambda^{+}}
  + \sum_{\mu\in N(\lambda)} \spinchar{\mu}^{*}.
 \]
 An analogous formula holds for the associate character
 $\spinchar{\lambda}^{a}$. In this case the character
 $\spinchar{\lambda^{+}}$ is replaced by its associate.
\end{theorem}
As stated in the previous section, in Section~\ref{sec:symm} it
will be important to consider the spin characters of the subgroups
$\theta^{-1}(\Scal_l\times\Scal_{n-l})$ for $1\le l \le n-1$. By
\cite[Theorem~4.3]{Stembridge89} every irreducible spin character
of $\theta^{-1}(\Scal_l\times\Scal_{n-l})$ is a \emph{reduced
Clifford product} $\spinchar{\mu}\otimes_z\spinchar{\nu}$ (see
\cite[Section~4]{Stembridge89} for a definition) for two strict
partitions $\mu\in\Dcal_l$ and $\nu\in\Dcal_{n-l}$. The character
values of $\spinchar{\mu}\otimes_z\spinchar{\nu}$ are readily
determined from the values of $\spinchar{\mu} \in \Irr_{-}(\tScal_l)$
and $\spinchar{\nu} \in \Irr_{-}(\tScal_{n-l})$ (see, for example,
\cite[Table~5.7]{HoffmanHumphreys} for a tabulation).

The \emph{projective outer product} is then simply defined to be the
induced character
\begin{equation*}
 \spinchar{\mu}\hat{\otimes}\spinchar{\nu}:=
 (\spinchar{\mu}\otimes_z\spinchar{\nu})\upa{\tScal_n}.
\end{equation*}
To decide which projective outer products are irreducible, we need
information on their constituents. Thanks to the work of Stembridge
in \cite{Stembridge89} these may be determined by an analogue of the
Littlewood-Richardson Rule (see \cite[2.8.13]{JamesKerber}) for spin
characters. In order to state this rule, we have to expand a little:

To a partition $\lambda\in\Dcal_n$ we associate a \emph{shifted diagram}
\begin{equation*}
 S(\lambda)=\{(i,j)\:|\: 1\le i\le \ell(\lambda),\: i\le j\le
 \lambda_i+i-1\}
\end{equation*}
whose elements we interpret as coordinates in a matrix style notation.
Let $A'=\{1'<1<2'<2<\dots\}$ be an ordered alphabet. The letters
$1',2',3',\dots$ are said to be \emph{marked}, the others are
\emph{unmarked}. We write $A$ for the set of unmarked letters of
$A'$. A \emph{shifted tableau} $T$ of shape $\lambda$ is a map $T:
S(\lambda)\rightarrow A'$ satisfying
\begin{itemize}
 \item[(1)] $T(i,j)\le T(i+1,j)$ for all $i,j$ with $(i,j),(i+1,j)\in
  S(\lambda)$; \hfill{(nondecreasing columns)}
 \item[(2)] $T(i,j)\le T(i,j+1)$ for all $i,j$ with $(i,j),(i,j+1)\in
  S(\lambda)$; \hfill{(nondecreasing rows)}
 \item[(3)] every $k\in\{1,2,\dots\}$ appears at most once in each 
  column of $T$;
 \item[(4)] every $k'\in\{1',2',\dots\}$ appears at most once in each row.
\end{itemize}
The \emph{content} of $T$ is a sequence of integers $(c_1,c_2,\dots)$,
where $c_k$ counts the number of nodes $(i,j)\in S(\lambda)$ such that
$T(i,j)$ is $k$ or $k'$.

For two strict partitions $\lambda$ and $\mu$ with $S(\mu)\subseteq
S(\lambda)$ the \emph{skew shifted diagram} $S(\lambda/\mu)$ of
shape $\lambda/\mu$ is the set $S(\lambda)\setminus S(\mu)$. The
corresponding \emph{skew shifted tableau} of shape $\lambda/\mu$ is
given by restricting the shifted tableau of shape $\lambda$ to this set.
Reading the rows of a (skew) shifted Tableau $T$ from left to right and
from bottom to top, gives its \emph{tableau-word} over $A'$.

We now count the number of skew shifted tableaux of shape $\lambda/\mu$
with content $\nu$ whose words are such that for all $a\in A$ the
leftmost letter of $\{a',a\}$ in $w$ is unmarked, and which fulfill
a set of further combinatorial conditions (called the \emph{lattice
property}) as detailed in \cite[Section~8]{Stembridge89}. For brevity,
we say that such a tableau satisfies (TP).

Denoting the number of tableaux of shape $\lambda/\mu$ with content
$\nu$ satisfying (TP) by $\mathrm{st}(\lambda;\mu,\nu)$, and setting
\begin{equation*}
 \varepsilon_{\alpha}:=
 \begin{cases}
  1, & \text{if } \alpha \text{ is even},\\
  \sqrt{2}, & \text{if } \alpha \text{ is odd}.
 \end{cases}
\end{equation*}
for a partition $\alpha$, allows us to finally state the
Littlewood-Richardson Rule.
\begin{theorem}[{\cite[Theorems 8.1, 8.3]{Stembridge89}, 
 \cite[Corollary 14.4]{HoffmanHumphreys}}]\label{thm:littlewoodrichardson}
 Let $\mu\in \Dcal_l$, $\nu\in \Dcal_{n-l}$ and $\lambda\in\Dcal_n$.
 Then we have
 \begin{equation}\label{sklprod}\tag{$\ast$} 
  \left(\spinchar{\mu} \hat\otimes \spinchar{\nu}, \spinchar{\lambda}
  \right) = \frac{1}{\varepsilon_{\lambda}\varepsilon_{\mu\cup\nu}}
  2^{\frac{1}{2}(\ell(\mu)+\ell(\nu)-\ell(\lambda))}\, \mathrm{st}(\lambda;
  \mu,\nu),
 \end{equation}
 unless $\lambda$ is odd and equal to $\mu\cup\nu$. In this case, 
 either $\spinchar\lambda$ or $\spinchar\lambda^a$ is constituent of
 $\spinchar\mu\hat\otimes\spinchar\nu$ with multiplicity 1. If
 $\ell(\lambda) > \ell(\mu)+\ell(\nu)$ then $\mathrm{st}(\lambda;\mu,\nu)$
 is zero, and $\spinchar{\lambda}$ is not a constituent.
\end{theorem}
\section{The Symmetric Groups}\label{sec:symm}
In this section we determine the irreducible spin characters of
$\tScal_n$ which are induced characters of proper subgroups.

In Section~\ref{sec:reductions} we gathered that there are
two types of subgroups left to consider: On the one hand, by
Lemma~\ref{lem:intransitive}, there are the preimages of the maximal
parabolic subgroups of $\Scal_n$, i.e., subgroups of the form
$\tScal_{l,n-l}:=\theta^{-1}(\Scal_l \times \Scal_{n-l})$. On the other
hand there are the subgroups $\tScal_{n:2}:=\theta^{-1}(\Scal_{n/2}\wr
\Scal_2)$ of Lemma~\ref{lem:transitive_imprimitive}.

\subsection{Characters Induced from $\tScal_{l,n-l}$}
Let us begin by considering spin characters induced from subgroups
of the first type. In other words, given two irreducible spin
characters $\spinchar{\mu} \in \tScal_l$ and $\spinchar{\nu} \in
\tScal_{n-l}$, we have to determine if their projective outer product
$\spinchar{\mu}\hat{\otimes}\spinchar{\nu}$ is again irreducible.
To decide this question, we make use of a classification of the
multiplicity-free projective outer products of two spin characters which
is mainly due to Bessenrodt.

To this end we need the following notation: A \emph{staircase}
is a partition of the form $(k,k-1,\dots,2,1)$, for some
$k\in\NN$. A \emph{fat staircase} is a partition of the form
$(k+r,k-1+r,\dots,2+r,1+r)$ for some $k\in\NN$, $r\ge 0$. In particular,
a staircase is also a fat staircase. A \emph{hook staircase} is the
concatenation of a fat staircase and a staircase (where one of them may
be empty).

In \cite[Theorem~3.2]{Bessenrodt02} Bessenrodt determined almost all
multiplicity-free projective outer products of $\tScal_n$. Together with
the first author's diploma thesis, this leads to the following theorem.
\begin{theorem}
\label{thm:mult-free}
 Let $\mu$ and $\nu$ be strict partitions of $l$ and $n-l$, respectively. The
 projective outer product 
 $\spinchar{\mu}\hat{\otimes}\spinchar{\nu}$ is multiplicity-free 
 if and only if it is as in one of the following cases:
 \begin{itemize}
  \item[(i)] $\spinchar{\text{\rm something}}\hat{\otimes}\spinchar1$.
  \item[(ii)] $\spinchar{\text{\rm hook~staircase}}\hat{\otimes}
    \spinchar{2,1}$, and the hook staircase is in $\Dcal_{n-3}^+$.
  \item[(iii)] $\spinchar{\text{\rm staircase}}\hat{\otimes}\spinchar m$ 
   for $m\in\NN$.
  \item[(iv)] $\spinchar{\text{\rm staircase}}\hat{\otimes}\spinchar{m-1,1}$ 
   for $m\in\NN$, $m>3$, such that the staircase and $(m-1,1)$ have 
   different signs.
  \item[(v)] $\spinchar{\text{\rm staircase}}\hat{\otimes}\spinchar{m+1,m}$ 
    for $m\in\NN$, and the staircase is in $\Dcal_{n-(2m+1)}^+$.
  \item[(vi)] $\spinchar{\text{\rm fat staircase}}\hat{\otimes}\spinchar m$ 
   for $m\in\NN$, $m>1$, such
   that the fat staircase and $(m)$ have different signs.
 \end{itemize}
\end{theorem}
While the proofs of cases (i) to (v) of Theorem~\ref{thm:mult-free} were
established in \cite{Bessenrodt02}, the proof of case (vi), which we
present here, was first conceived in \cite{Nett07}.
\begin{proof}
 Let $\mu=(\mu_1,\dots,\mu_{\ell(\mu)})$ and
 $\nu=(\nu_1,\dots,\nu_{\ell(\nu)})$ be strict partitions, and assume
 without loss of generality that $\ell(\mu)\ge \ell(\nu)$. Then the partition
 \[
  \lambda=\mu+\nu=(\mu_1+\nu_1,\dots,\mu_{\ell(\nu)}+\nu_{\ell(\nu)},
  \mu_{\ell(\nu)+1},\dots,\mu_{\ell(\mu)})
 \] 
 yields a constituent of $\spinchar{\mu}\hat{\otimes}\spinchar{\nu}$,
 and the factor $2^{\frac{1}{2}(\ell(\mu)+\ell(\nu)-\ell(\lambda))}$
 in formula \eqref{sklprod} of the Littlewood-Richardson
 rule simplifies to $2^{\frac{1}{2}\ell(\nu)}$. Now, assuming that
 $\spinchar{\mu}\hat{\otimes}\spinchar{\nu}$ is multiplicity-free, and
 since there is exactly one tableau of shape $\lambda/\mu$ with content
 $\nu$, the factor $2^{\frac{1}{2}\ell(\nu)}$ is at most two, and
 therefore $1\le\ell(\nu)\le 2$.

 Also, as $\spinchar{\mu}\hat\otimes\spinchar{\nu}$ is multiplicity free, the
 inequality $\mathrm{st}(\lambda;\mu,\nu) \le 2$ holds for all $\lambda \in
 \Dcal_n$. In particular, if we have $\mathrm{st}(\lambda;\mu,\nu)\le 1$ 
 for all $\lambda \in \Dcal_n$, then one of the cases (i)-(v) of 
 Theorem~\ref{thm:mult-free}
 applies, 
 and we are done by \cite[Theorem~3.2]{Bessenrodt02}.
 Thus we consider the previously untreated case that there exists a
 $\lambda_0 \in \Dcal_n$ such that $\mathrm{st}(\lambda_0;\mu,\nu)=2$.
 This implies that the coefficient of \eqref{sklprod} simplifies
 to $\frac{1}{\varepsilon_{\lambda_0} \varepsilon_{\mu\cup\nu}}
 2^{\frac{1}{2}(\ell(\mu)+\ell(\nu)-\ell(\lambda_0))}=\frac{1}{2}$,
 and hence $\ell(\lambda_0)=\ell(\mu)+\ell(\nu)$ and
 $\varepsilon_{\lambda_0}=\varepsilon_{\mu\cup\nu}=\sqrt{2}$, i.e.,
 the partitions $\mu$ and $\nu$ have different signs.

 Let us assume $\ell(\nu)=2$. If $\mu$ is a staircase
 or $\nu=(2,1)$, there is no partition $\lambda$ with
 $\mathrm{st}(\lambda;\mu,\nu)>0$ and $\ell(\lambda)= \ell(\mu)+2$. On the
 other hand, if $\mu$ is not a staircase and $\nu\neq(2,1)$, there exists a
 partition $\lambda$ with $\mathrm{st}(\lambda;\mu,\nu)\ge 2$ and
 $\ell(\lambda)<\ell(\mu)+2$, giving 
 $(\spinchar{\lambda},\spinchar{\mu}\hat\otimes\spinchar{\nu})>1$.
 
 Therefore we conclude $\ell(\nu)=1$, and we may assume that we are not in
 one of the Cases (i) to (v). In particular, we can assume that
 $\mu$ is not a staircase. Now, if $\mu_{i+1}+1<\mu_i$ for some
 $i<\ell(\mu)$, there exists at least one partition $\lambda$ with
 $\mathrm{st}(\lambda;\mu,\nu)\ge 2$ and $\ell(\lambda)<\ell(\mu)+1$, again 
 giving rise to constituent whose multiplicity exceeds one.
 Thus $\mu$ has to be a fat staircase, and we are in case (vi).

 For the converse, it is easy to see that in Cases
 (i) to (vi) the projective outer product
 $\spinchar\mu\hat\otimes\spinchar\nu$ is multiplicity-free.
\end{proof}

As a corollary to the Branching Rule~\ref{thm:BR} and the previous
theorem, we can now determine the irreducible projective outer products
$\spinchar\mu\hat\otimes\spinchar\nu$.

\begin{lemma}\label{la:parabolSn}
 Let $\mu\in\Dcal_l$ and $\nu\in\Dcal_{n-l}$. 
 \begin{itemize}
  \item[(a)] If $l=n-1$, then the spin character
   $\spinchar\mu\hat\otimes\spinchar1= \spinchar{\mu}\upa{\tScal_n}$
   is irreducible if and only if $\mu=(m,m-1,\dots,1)$ for some
   $m\in\NN$ with $m\equiv 2$ or $3\pmod 4$. In this case we
   have $\mu\in\Dcal_{n-1}^-$ and $\spinchar\mu\upa{\tScal_n}
   =\spinchar\lambda$ with $\lambda=(m+1,m-1,m-2,\dots,1)\in\Dcal_n^+$.
  \item[(b)] If $n-l$ and $l$ are larger than $1$, then
   $\spinchar{\mu}\hat{\otimes}\spinchar{\nu}$ is a
   reducible character of $\tScal_n$. In particular,
   $\spinchar{\mu}\hat{\otimes}\spinchar{\nu}$ contains
   (at least) two non-associate irreducible constituents,
   except when the partition $\mu=(m,m-1,\ldots,1)$ is an
   even staircase and $\nu=(2,1)$. In this case
   $\spinchar{\mu}\hat{\otimes}\spinchar{\nu}=\spinchar{\lambda}+
   \spinchar{\lambda}^a$ with $\lambda=(m+2,m,m-2,\dots,2,1)$.
 \end{itemize}
\end{lemma}
\begin{proof}    
  (a) 
 Let $\spinchar{\mu}$ be a spin character of $\tScal_{n-1}$
 such that $\spinchar{\mu}\upa{\tScal_n}=\spinchar{\lambda}$ is
 irreducible. Then in the terminology of Theorem~\ref{thm:BR} we have
 $|N(\mu)|=1$ and $\mu_{\ell(\mu)}=1$, i.e., $\mu=(m,m-1,\dots,1)$
 for $m=\ell(\mu)$. Since $\spinchar{\lambda}^*$ appears as a summand
 in $\spinchar{\mu}\upa{\tScal_n}$ by the Branching Rule, we have
 $\lambda\in\Dcal_n^+$ and hence $\mu \in \Dcal_{n-1}^-$. It follows
 that $\lambda=(m+1,m-1,m-2,\dots,1)$. Since $n=1+\sum_{i=1}^m
 i=m(m+1)/2+1$ and $n-m=m(m-1)/2+1$ is even, $m$ is congruent to $2$ or
 $3$ modulo $4$.

  (b)
 We consider the Cases (ii)-(vi) of Theorem \ref{thm:mult-free}. 

 Case (ii):
 Let $\mu\in \Dcal_{n-3}^+$ be a hook staircase and $\nu=(2,1)$.
 We write $\mu=\mu^{(1)}\sqcup \mu^{(2)}$ for a fat
 staircase $\mu^{(1)}=(k+r,k-1+r,\dots,1+r)$ and a staircase
 $\mu^{(2)}=(m,m-1,\dots,1)$. Note that one of these parts may be empty.
 If $\mu^{(1)}$ is not empty then $r\ge 1$; If both $\mu^{(1)}$ and
 $\mu^{(2)}$ are non-empty then $\mu^{(1)}_k=1+r>m+1 = \mu^{(2)}_1+1$.
 
 Assume that $\mu^{(1)}$ is not empty and let $k\ge
 2$, then the spin character $\spinchar{\lambda}$ is a constituent of 
 $\spinchar{\mu}\hat{\otimes}\spinchar{2,1}$ for 
 \begin{itemize}
  \item[(a)] $\lambda=(\mu^{(1)}_1+2,\mu^{(1)}_2+1,\mu^{(1)}_3,\dots)$, and
  \item[(b)] $\lambda=(\mu^{(1)}_1+2,\mu^{(1)}_2,\mu^{(1)}_3,
   \dots, \mu^{(1)}_k, \mu^{(2)}_1+1, \mu^{(2)}_2,
   \mu^{(2)}_3,\ldots,\mu^{(2)}_m)$
   
   (with $\lambda=
   (\mu^{(1)}_1+2,\mu^{(1)}_2+1,\mu^{(1)}_3,\dots,\mu^{(1)}_k, 1)$ if
   $\mu^{(2)}$ is empty).
 \end{itemize}

 Let $k=1$, i.e., $\mu^{(1)}=(r+1)$. If $\mu^{(2)}$ is not empty, then
 the spin character $\spinchar{\lambda}$ for $\lambda$ as in (b) above
 is again a constituent of $\spinchar{\mu}\hat{\otimes}\spinchar{2,1}$.
 Furthermore, $\spinchar{\lambda}$ with
 $\lambda=(\mu^{(1)}_1+1,\mu^{(2)}_1+1,\mu^{(2)}_1,\mu^{(2)}_3,
 \dots,\mu^{(2)}_m)$ is a constituent, too. If $\mu^{(2)}$
 is empty, the characters $\spinchar{\mu^{(1)}_1+2,1}$
 and $\spinchar{\mu^{(1)}_1+1,2}$ are constituents of
 $\spinchar{\mu}\hat{\otimes}\spinchar{2,1}$.

 Now assume that $\mu^{(1)}$ is empty. Let
 $\mu=\mu^{(2)}=(m,m-1,\dots,1)$ be an even staircase. Hence $m$
 is congruent to $0$ or $1$ modulo $4$, and therefore $m \ge
 4$. In this case only $\lambda=(m+2,m,m-2,\dots,1)$ yields a
 constituent of the outer projective product. The corresponding
 spin character $\spinchar{\lambda}$ is not self-associate and
 therefore we have $\spinchar{\mu}\hat{\otimes}\spinchar{2,1}=
 \spinchar{\lambda}+\spinchar{\lambda}^a$.
 
 Case (iii):
 Let $\mu=(k,k-1,\dots,1)$ be a staircase and $\nu=(m)$ for $k$
 and $m>1$. In this case, the partitions $(k+m,k-1,\dots,1)$ and
 $(k+m-1,k,k-2,\dots,1)$ yield two non-associate constituents of
 $\spinchar{\mu}\hat{\otimes}\spinchar{m}$.

 Case (iv):
 If $\mu=(k,k-1,\dots,1)$ is again a staircase with $k>1$
 and $\nu=(m-1,1)$ for $m>3$, then the two partitions
 $(k+m-1,k,k-2,\dots,1)$ and $(k+m-2,k+1,k-2,\dots,1)$ yield
 constituents of $\spinchar{\mu}\hat{\otimes}\spinchar{m-1,1}$.

 Case (v):
 Let $\mu=(k,k-1,\dots,1)$ be an even staircase and
 $\nu=(m+1,m)$ for $m>1$ (we have already considered the case
 $\nu=(2,1)$ and $\mu$ a staircase). By the same argument as in Case
 (ii) above, we have $k\equiv 0$ or $1 \pmod 4$. Since $k>1$, this
 implies $k\ge 4$. The partitions $(k+m+1,(k-1)+m,k-2,\dots,1)$ and
 $(k+m+1,k+m-2,k-1,k-3,\dots,1)$ yield constituents of the projective
 outer product $\spinchar{\mu}\hat{\otimes}\spinchar{m+1,m}$.

 Case (vi):
 Let $\mu=(k+r, k-1+r,\dots,1+r)$ be a fat staircase for some $k,r\ge
 1$ and $\nu=(m)$ for some $m>1$. Then the partitions $(k+r+m,
 k-1+r,\dots,1+r)$ and $(k+r+m-1,k-1+r,\dots,1+r,1)$ yield two
 non-associate constituents of $\spinchar\mu\hat\otimes\spinchar\nu$.
\end{proof}

\subsection{Characters Induced from $\tScal_{n:2}$}
We will now consider the imprimitive and transitive subgroups
$\tScal_{n:2}$ of $\tScal_n$. Let $n = 2m$, where $m\ge 2$. The subgroup
$\tScal_{m,m}$ has index two in $\tScal_{n:2}$, and therefore the
irreducible characters of $\tScal_{n:2}$ may be determined through an
elementary application of Clifford Theory. Every irreducible character
of $\tScal_{m,m}$ gives rise to one or two irreducible characters of
$\tScal_{n:2}$, and all irreducible characters of $\tScal_{n:2}$ arise
in this manner: A character is either invariant under the conjugation
action of $\tScal_{n:2}$ and thus possesses two distinct extensions
to $\tScal_{n:2}$, or its inertia subgroup is $\tScal_{m,m}$ and
the character fuses with its conjugate to give a single irreducible
character of $\tScal_{n:2}$ by induction.

For $\spinchar{\mu}, \spinchar{\nu} \in \Irr_{-}(\tScal_m)$
let $(\spinchar{\mu}\otimes_z \spinchar{\nu})^+$ and
$(\spinchar{\mu}\otimes_z \spinchar{\nu})^-$ denote the two distinct
extensions of the reduced Clifford product $\spinchar{\mu}\otimes_z
\spinchar{\nu}$, if the latter is invariant. Analogously, we denote by
$(\spinchar{\mu}\otimes_z \spinchar{\nu})^0$ the irreducible induced
character $(\spinchar{\mu}\otimes_z \spinchar{\nu})\upa{\tScal_{n:2}}$,
if $\spinchar{\mu}\otimes_z \spinchar{\nu}$ is not invariant.

The following first result is immediate.
\begin{corollary}
 \label{cor:notinvariant}
 Let $(\spinchar{\mu}\otimes_z\spinchar{\nu})^0 \in \Irr_{-}(\tScal_{n:2})$.
 Then the induced character
 $(\spinchar{\mu}\otimes_z\spinchar{\nu})^0\upa{\tScal_n}$ is always
 reducible.
\end{corollary}
\begin{proof}
 This is a simple consequence of Lemma~\ref{la:parabolSn} and the
 transitivity of induction.
\end{proof}
The analysis if either of the characters $(\spinchar{\mu}\otimes_z
\spinchar{\nu})^+$ and $(\spinchar{\mu}\otimes_z \spinchar{\nu})^-$
induces irreducibly to $\tScal_n$ is slightly more involved. As both are
extensions of an invariant character of $\tScal_{m,m}$ we will determine
these first.

Let $\Scalm := \left< s_1,\ldots,s_{m-1}\right>$ and
$\Scalmm:=\left<s_{m+1},\ldots,s_{2m-1}\right>$. The image of
$\tScal_{n:2}$ under $\theta$ is isomorphic to $(\Scalm \times \Scalmm)
\rtimes \left< \tau \right>$, for which the action of $\tau$ induces
the outer automorphism which maps $s_i$ to $s_{i+m}$, if the indices
are taken modulo $m$. Let $\ttau$ denote the standard lift of $\tau$
in $\tScal_n$. Then a character of $\tScal_{m,m}$ is invariant under
the action of $\tScal_{n:2}$ if and only if it is invariant under
conjugation by $\ttau$. Therefore we consider the action of $\ttau$ on
the conjugacy classes of $\tScal_{m,m}$.

As $t_it_j = zt_jt_i$ for $|i-j|>1$, any element $y$ of $\tScal_{m,m}$
may be written as a product $gh$ for some $g \in \tScalm
:= \left< t_1,\ldots,t_{m-1}\right>$ and $h\in \tScalmm :=
\left< t_{m+1},\ldots,t_{2m-1}\right>$, and we have $gh=
z^{\sigma(g)\sigma(h)}hg$. If the two elements $y$ and $zy$
are not conjugate in $\tScal_{m,m}$, the full preimage of
$\theta(y)^{\Scal_m\times \Scal_m}$ is the union of the two classes
containing $y$ and $zy$. In this case we say that the conjugacy class
of $\theta(y)$ \emph{splits}, or that the classes of $y$ and $zy$ are
\emph{split}. As a class of $\Scal_m \times \Scal_m$ is naturally
parameterized by a pair $(\pi, \mu)$ of partitions of $m$, the class
$y^{\tScal_{m,m}}$ is therefore also parameterized by $(\pi, \mu)$, if
the cycle types of $\theta(g)$ and $\theta(h)$ are $\pi$ and $\mu$,
respectively. By a slight abuse of notation we denote conjugacy classes
by their parameters. If there are two classes with the same parameter,
as is the case for split classes, we affix subscripts to distinguish
them.

We begin by determining the classes of $\Scal_m \times \Scal_m$ which
split. Let $\Ocal_m$ denote the set of all odd part partitions of $m$.
\begin{lemma}
 \label{la:splittingclasses}
 A class $(\pi, \mu)$ of $\Scal_m\times \Scal_m$ splits if and only if
 $(\pi, \mu)$ is an element of $\Ocal_m \times \Ocal_m$, $\Dcal^-_m
 \times \Dcal_m^+$, or $\Dcal_m^+\times \Dcal_m^-$.
\end{lemma}
\begin{proof}
 If $(\pi,\mu) \in \Ocal_m \times \Ocal_m$, then $\theta(gh)$ has odd
 order. Hence we may assume that $gh$ has odd order, too. For any $y
 \in \theta^{-1}(C_{\Scal_m\times \Scal_m}(\theta(gh)))$ we have that
 $(gh)^y \in \{ gh, zgh\}$, so as the order of $gh$ is odd, we conclude
 $(gh)^y =gh$. Therefore the centralizer of $gh$ in $\tScal_{m,m}$ is
 the full preimage of $C_{\Scal_m\times \Scal_m}(\theta(gh))$ and the
 class of $\theta(gh)$ splits.

 If $(\pi,\mu) \in \Dcal_m^- \times \Dcal_m^+$, then $\pi$
 possesses an odd number $o$ of even parts, and $\mu$ has an
 even number $e$ of even parts. For $i=1,\ldots,\ell(\pi)$
 let $\tilde{C}_{\pi_i}\in \tScalm$ denote a lift of a
 $\pi_i$-cycle. Likewise, for $j= 1,\ldots, \ell(\mu)$ let
 $\tilde{C}_{\mu_i} \in \tScalmm$ be a lift of a $\mu_j$-cycle, and
 set $g:=\tilde{C}_{\pi_1}\cdots\tilde{C}_{\pi_{\ell(\pi)}}$ and
 $h:=\tilde{C}_{\mu_1}\cdots\tilde{C}_{\mu_{\ell(\mu)}}$. If $\pi_i$ is
 odd, then $\tilde{C}_{\pi_i}gh = gh\tilde{C}_{\pi_i}$, and if $\pi_i$
 is even, we have $\tilde{C}_{\pi_i}gh = z^{o+e-1}gh\tilde{C}_{\pi_i}
 = gh\tilde{C}_{\pi_i}$, as $o+e-1$ is even. Thus $\tilde{C}_{\pi_i}
 \in C_{\tScal_{m,m}}(gh)$ for all $i=1,\ldots,\ell(\pi)$. The same
 holds for all $\tilde{C}_{\mu_j}$, $j=1,\ldots,\ell(\mu)$. Therefore
 $C_{\tScal_{m,m}}(gh)$ is again the full preimage of the centralizer of
 $\theta(gh)$, and thus the class of $\theta(gh)$ splits.

 Interchanging the roles of $\pi$ and $\mu$ above, yields
 the result for $(\pi,\mu) \in \Dcal_m^+\times\Dcal_m^-$. By
 \cite[Theorem~5.9]{HoffmanHumphreys} there are exactly $|\Dcal_m|^2 +
 2|\Dcal_m^+||\Dcal_m^-|$ splitting classes. As $|\Ocal_m| = |\Dcal_m|$
 we are done.
\end{proof}
In order to give the action of $\ttau$ on the conjugacy classes
of $\tScal_{m,m}$, we first study its effect on the generators
$t_1,\ldots,t_{m-1},t_{m+1},\ldots,t_{2m-1}$ of $\tScal_{m,m}$. Note
that the group $\tScal_{m,m}$ also possesses the outer automorphism
$\iota$ which maps $t_i$ to $t_{i+m}$ where the indices are again taken
modulo $m$. There is a subtle difference between the actions of $\ttau$
and $\iota$ depending on the parity of $m$.

\begin{lemma}\label{la:ttauonti}
 Let $1\le i \le m-1$. Then $t_i^{\ttau} = z^m t_{i+m}$ and
 $t_{i+m}^{\ttau} = z^m t_{i}$.
\end{lemma}
\begin{proof}
 For $1\le k \le m$ we set
 \[
  \tilde{\tau}_k := t_{k+m-1} t_{k+m-2}\cdots t_{k+1}
  t_k t_{k+1}\cdots t_{k+m-2}t_{k+m-1}. 
 \]
 Then $\ttau = \ttau_1\cdots \ttau_m$. First note that we have
 $\ttau_i\ttau_j = z\ttau_j\ttau_i$ for all $1\le i,j\le m$: The element
 $\theta(\ttau_i\ttau_j)$ of $\Scal_n$ has cycle type $(2^2,1^{2m-4})
 \notin \Ocal_n \cup \Dcal_n^-$. Thus $\ttau_i\ttau_j$ is conjugate to
 $t_1t_3$ and therefore has order four. This yields $\ttau_i\ttau_j
 = z\ttau_j^{-1}\ttau_i^{-1}$, which is $z\ttau_j\ttau_i$ as every
 $\ttau_k$ is the product of $2m-1$ factors. Using the relations given
 in the presentation of $\tScal_n$ we obtain $\ttau^{-1} t_i \ttau = z^m
 \ttau_{i+1}\ttau_{i} t_i \ttau_{i}\ttau_{i+1}$ for $1\le i \le m-1$.
 Further meticulous applications of the these relations give that the
 latter is equal to $\ttau^{-1}t_i\ttau = z^m t_{i+m}$, as claimed.
 Therefore the equation $\ttau^{-1}t_{i+m}\ttau = z^m\ttau^{-2} t_i
 \ttau^2 = z^m t_i$ follows, too.
\end{proof}
\begin{lemma}\label{la:iota}
 Let $\iota$ denote the automorphism of $\tScal_{m,m}$ taking $t_i$
 to $t_{i+m}$ (where the indices are taken modulo $2m$). In other
 words, for $g \in \tScalm$ and $h \in \tScalmm$ we have $(gh)^\iota
 = z^{\sigma(g)\sigma(h)} h^\iota g^\iota$. Hence $\iota$ maps the
 conjugacy class containing the standard lift of $\theta(gh)$ to the
 conjugacy class containing the standard lift of $\theta(h^\iota
 g^\iota)$.
\end{lemma}
\begin{proof}
 For a non-split class the assertion is clear. Let $gh$ lie in a split
 class, and without loss of generality assume that $gh$ is the standard
 lift of $\theta(gh)$, i.e., either both $g$ and $h$ are the standard
 lifts of $\theta(g)$ and $\theta(h)$, or both are not. The same
 holds for the image under $\iota$: By the definition of $\iota$ both
 $g^\iota$ and $h^\iota$ are standard lifts if and only if both $g$
 and $h$ are standard lifts. Now as the element $\theta(gh)$ lies in a
 conjugacy class of $\Scal_m\times \Scal_m$ with parameter $(\pi,\mu)$,
 where $(\pi,\mu)$ is an element of either $\Ocal_m \times \Ocal_m$,
 $\Dcal_m^-\times \Dcal_m^+$, or $\Dcal_m^+ \cap \Dcal_m^-$, we have
 that $\sigma(g) = 0$ or $\sigma(h) = 0$. Therefore $(gh)^\iota$ is the
 standard lift.
\end{proof}

\begin{corollary}\label{cor:actionttauonclasses}
 For a splitting class of $\Scal_m\times \Scal_m$ denoted by its
 parameter $(\pi,\mu)$ let $(\pi,\mu)_1$ be the class of $\tScal_{m,m}$
 containing the standard lift $gh$ of the canonical representative of
 $(\pi,\mu)$. Accordingly, let $(\pi,\mu)_2$ denote the class containing
 $zgh$.\\[1mm]
 If $m$ is even, then $(\pi,\mu)_1^{\ttau} = (\mu,\pi)_1$.\\[1mm]
 If $m$ is odd, 
 \[ (\pi,\mu)_1^{\ttau} = 
  \begin{cases}
   (\mu,\pi)_1, & \text{if }(\pi,\mu) \in \Ocal_m \times \Ocal_m,\\
   (\mu,\pi)_2, & \text{if }(\pi, \mu) \in \Dcal_m^- \times \Dcal_m^+
   \text{ or }\Dcal_m^+\times \Dcal_m^-.
  \end{cases}
 \]
\end{corollary}
\begin{proof}
 If $m$ is even, then by Lemma~\ref{la:ttauonti} we see that $\ttau$ and
 $\iota$ act identically, and Lemma~\ref{la:iota} gives the claim. If
 $m$ is odd, by Lemma~\ref{la:ttauonti} we have
 \[ 
  (gh)^{\ttau} = z^{\sigma(g)+\sigma(g)\sigma(h)+\sigma(h)} h^\iota
  g^\iota.
 \] 
 Thus if $gh \in (\pi,\mu)_1$ for $(\pi,\mu) \in \Dcal_m^- \times
 \Dcal_m^+$ or $\Dcal_m^+ \times \Dcal_m^-$, then $(gh)^{\ttau} =
 zh^\iota g^\iota$. Hence by Lemma~\ref{la:iota} we have $(gh)^{\ttau} =
 z(gh)^\iota\in (\mu,\pi)_2$.
\end{proof}
Summing up, we can give the action of $\ttau$ on the irreducible
characters of $\tScal_{m,m}$ as follows:
\begin{corollary}
 \label{cor:actionttauonchars}
 For $\chi \in \Irr_{-}(\tScal_{m,m})$ the action of $\ttau$ is given by
 \[
  \chi^{\ttau}=
  \begin{cases}
   \chi^\iota, & \text{if } $m$ \text{ is even},\\
   (\chi^\iota)^a, & \text{if } $m$ \text{ is odd}.
  \end{cases}
 \]
\end{corollary}
\begin{proof}
 This is now immediate by Corollary~\ref{cor:actionttauonclasses}.
\end{proof}
With the action of $\ttau$ on $\Irr_{-}(\tScal_{m,m})$ known, it is easy
to determine which characters are $\ttau$-invariant.
\begin{corollary}\label{cor:invariant}
 Let $\mu, \nu \in \Dcal_m$. Then
 $\spinchar{\mu}\otimes_z\spinchar{\nu}$ is $\ttau$-invariant if and
 only if, $\mu = \nu$.
\end{corollary}
\begin{proof}
 Let $\chi:= \spinchar{\mu}\otimes_z \spinchar{\nu}$, then by
 Corollary~\ref{cor:actionttauonchars} we have $\chi^{\ttau}\downa{\ker
 \sigma} = \chi^{\iota}\downa{\ker \sigma}$ independent of the
 parity of $m$. If $\chi$ is $\ttau$-invariant, this implies
 $(\spinchar{\mu}\otimes_z\spinchar{\nu})\downa{\ker \sigma} =
 (\spinchar{\nu} \otimes_z \spinchar{\mu})\downa{\ker \sigma}$.
 We may now argue as in \cite[proof of 5.9]{HoffmanHumphreys}:
 From the character values of the reduced Clifford product
 (cf.\ \cite[Table~5.7]{HoffmanHumphreys}) it follows that
 $\spinchar{\mu}\downa{\ker \sigma} = \spinchar{\nu}\downa{\ker
 \sigma}$. Hence $\spinchar{\mu} = \spinchar{\nu}$ or $\spinchar{\mu} =
 \spinchar{\nu}^a$, so $\mu = \nu$.

 The converse follows from the fact that for any spin character
 $\spinchar{\mu} \in \Irr_{-}(\tScal_m)$ the reduced Clifford product
 $\spinchar{\mu}\otimes_z \spinchar{\mu}$ is self-associate. Hence it is
 $\ttau$-invariant by Corollary~\ref{cor:actionttauonchars}.
\end{proof}
Having the $\tScal_{n:2}$-invariant characters of $\tScal_{m,m}$ at our
disposal, we can study their behavior under induction.
\begin{lemma}\label{la:indvalues}
 Let $\chi \in \Irr_{-}(\tScal_{m,m})$ be invariant under the
 conjugation action of $\tScal_{n:2}$, and let $\chi\upa{\tScal_{n:2}}
 = \chi^+ + \chi^-$, where $\chi^+, \chi^- \in \Irr_{-}(\tScal_{n:2})$
 are two extensions of $\chi$. Then for all $g \in \tScal_n$ we
 have $\chi^+\upa{\tScal_n}(g) = \chi^-\upa{\tScal_n}(g)$, except
 possibly when $\theta(g)$ has cycle type $2\pi \in \Dcal_n^-$. In
 this case $\chi^-\upa{\tScal_n}(g) = \chi^+\upa{\tScal_n}(zg) =
 -\chi^+\upa{\tScal_n}(g)$.
\end{lemma}
\begin{proof}
 Let $C$ denote a conjugacy class of $\tScal_n$. The value of
 $\chi^+\upa{\tScal_n}$ on $C$ is determined by the values of $\chi^+$
 on the conjugacy classes of $\tScal_{n:2}$ lying in $C$. The inner and
 outer classes do not interfere in the following sense: If there are
 outer classes of $\tScal_{n:2}$ lying in $C$, then the character value
 of $\chi^+\upa{\tScal_n}$ on $C$ depends only on the values of $\chi^+$
 on these outer classes. To this end let $(\pi, \mu)$ denote an inner
 class of $\tScal_{n:2}$ fusing via induction to $\tScal_{n}$ into $C$.
 As by \cite[4.2.17]{JamesKerber} the outer classes are parameterized
 by $2\Pcal_{m}$, we conclude that both $\pi$ and $\mu$ only consist
 of even parts. By \cite[Table~5.7]{HoffmanHumphreys} $\chi^+$ may
 only assume nonzero values on $(\pi, \mu)$ if both partitions are
 even. Furthermore, setting $\chi =\spinchar{\lambda} \otimes_z
 \spinchar{\lambda}$, its value on $(\pi, \mu)$ is up to a constant the
 product of the values of the spin character $\spinchar{\lambda}\in
 \Irr_{-}(\tScal_m)$ on the corresponding classes of $\tScal_{m}$ with
 parameters $\pi$ and $\mu$. Since a spin character of $\tScal_m$ only
 takes nonzero values on classes parameterized by elements of $\Ocal_m$
 or $\Dcal^-_m$, we conclude that $\chi^+$ is zero on all inner classes
 which fuse into $C$.

 It is now immediate that $\chi^+\upa{\tScal_n}(g)
 =\chi^-\upa{\tScal_n}(g)$ for all $g \in \tScal_n$, except when
 $\theta(g)$ has cycle type $2\pi$, i.e., when an outer class of
 $\tScal_{n:2}$ lies in the $\tScal_n$-conjugacy class of $g$.
 Now, if $\theta(g)$ has cycle type $2\pi \notin \Dcal_n^-$, the
 elements $g$ and $zg$ are conjugate in $\tScal_n$, and therefore
 $\chi^+\upa{\tScal_n}(g) = 0 = \chi^-\upa{\tScal_n}(g)$.
 So let $\theta(g)$ have cycle type $2\pi \in \Dcal_n^-$
 for an element $g \in \tScal_{n:2}$. Then the two classes
 $g^{\tScal_{n:2}}$ and $(zg)^{\tScal_{n:2}}$ fuse into the two
 classes $g^{\tScal_n}$ and $(zg)^{\tScal_n}$, respectively.
 Setting $e := |C_{\tScal_n}(g):C_{\tScal_{n:2}}(g)|$, we arrive
 at $\chi^-\upa{\tScal_n}(g) = e\chi^-(g) = -e\chi^+(g) =
 \chi^+\upa{\tScal_n}(zg)$, as the extensions $\chi^+$ and $\chi^-$
 fulfill $\chi^-(g) = -\chi^+(g)$.
\end{proof}

In the light of Lemma~\ref{la:indvalues} it is not surprising
that either both induced extensions $(\spinchar{\mu}\otimes_z
\spinchar{\mu})^+$ and $(\spinchar{\mu}\otimes_z \spinchar{\mu})^-$ are
irreducible, or both are not.
\begin{corollary}\label{cor:indpairsSn}
 Under the hypothesis of Lemma~\ref{la:indvalues} the character
 $\chi^+\upa{\tScal_n}$ is irreducible if and only if
 $\chi^-\upa{\tScal_n}$ is irreducible.
\end{corollary}
\begin{proof}
 By Lemma~\ref{la:indvalues} the induced characters
 $\chi^+\upa{\tScal_n}$ and $\chi^-\upa{\tScal_n}$ have the same norm.
\end{proof}

\begin{lemma}
 \label{la:kleinevielfachheit}
 Let $\mu\in\Dcal_m$. The spin character
 $\spinchar{\mu}\hat{\otimes}\spinchar{\mu}$ contains (at
 least) two non-associate, distinct constituents or one
 constituent with multiplicity four, except when $m=3$ and
 $\mu=(2,1)$, or $m=6$ and $\mu=(3,2,1)$. In these cases we have
 $\spinchar{2,1}\hat\otimes\spinchar{2,1}=2\spinchar{4,2}$ and
 $\spinchar{3,2,1}\hat\otimes\spinchar{3,2,1}= 2(\spinchar{6,4,2}+
 \spinchar{6,4,2}^a)$.
\end{lemma} 
\begin{proof}
 Taking $\lambda = 2\mu$ yields a constituent $\spinchar{\lambda}$
 of $\spinchar \mu \hat\otimes \spinchar\mu$, and the factor
 $2^{\frac{1}{2}(2\ell(\mu)-\ell(\lambda))}$ in formula \eqref{sklprod}
 of \ref{thm:BR} simplifies to $2^{\frac{1}2 \ell(\mu)}$. Note
 that $\varepsilon_{\mu\cup\mu}=1$ since $\mu\cup\mu$ is an even
 partition. If $\ell(\mu)\ge 4$ then the coefficient $2^{\frac12
 \ell(\mu)}\frac1{\varepsilon_{\lambda}}$ in \eqref{sklprod} is at least
 4. If $\ell(\mu)=1$, i.e., $\mu=(m)$ for some $m\ge 2$, then $(2m-1,
 1)$ yields a second constituent not associate to $\spinchar{2m}$.
 If $\ell(\mu)=2$ then $\lambda=2\mu$ is an even partition and the
 factor 2 appears in formula \eqref{sklprod}. Hence, we may assume
 $\mathrm{st}(2\mu;\mu,\mu) = 1$. By \cite[Theorem~2.2]{Bessenrodt02}
 this leaves $\mu=(2,1)$, and the corresponding projective outer
 product is $2 \spinchar{4,2}$. If $\ell(\mu)=3$ an analogous argument
 gives $\mu=(3,2,1)$, and we have $\spinchar{3,2,1} \hat\otimes
 \spinchar{3,2,1} = 2(\spinchar{6,4,2}+ \spinchar{6,4,2}^a)$.
\end{proof}

\begin{corollary}
 \label{cor:inducedimprimitive}
 If $n=6$ and $\mu=(2,1)$ we have $(\spinchar\mu
 \otimes_z\spinchar\mu)^+\upa{\tScal_6} = (\spinchar\mu
 \otimes_z\spinchar\mu)^-\upa{\tScal_6} =\spinchar{4,2}
 \in\Irr_-(\tScal_6)$. In all other cases, both
 $(\spinchar{\mu}\otimes_z \spinchar{\mu})^+$ and
 $(\spinchar{\mu}\otimes_z \spinchar{\mu})^-$ of
 $\Irr_{-}(\tScal_{n:2})$ induce reducibly to $\tScal_n$
\end{corollary}
\begin{proof}
 We have $\spinchar{\mu}\hat{\otimes} \spinchar{\mu} =
 (\spinchar{\mu}\otimes_z \spinchar{\mu})^+\upa{\tScal_n} +
 (\spinchar{\mu}\otimes_z \spinchar{\mu})^-\upa{\tScal_n}$. By
 Corollary~\ref{cor:indpairsSn} either both extensions induce
 irreducibly, or both do not. The only possibility for $\mu$ such that
 $\spinchar{\mu}\hat\otimes\spinchar{\mu}$ has norm 2 is given by
 Lemma~\ref{la:kleinevielfachheit} to be $\mu = (2,1)$.
\end{proof} 
With the help of Lemma~\ref{la:kleidmanwales}, Lemma~\ref{la:parabolSn}
and Corollary~\ref{cor:inducedimprimitive} we may readily identify
all imprimitive faithful characters of $\tScal_n$, and state the
corresponding maximal subgroups which occur as their block stabilizers.
For convenience, we summarize our current standing in the following
corollary.
\begin{corollary}\label{cor:mainSn}
 Let $\spinchar{\lambda}$ be an imprimitive spin character of $\tScal_n$.
 Then
 \begin{itemize}
  \item[(i)] $\lambda \in \Dcal_n^+$, and $\spinchar{\lambda}$
   is the induced of either constituent of 
   $\spinchar{\lambda}\downa{\tAcal_n}$, i.e., $\tAcal_n$ is a block
   stabilizer.
 \end{itemize}
   Furthermore,
   \begin{itemize}
  \item[(ii)]  
   if $\lambda = (m+1,m-1,\ldots,1)\in \Dcal_n^+$ for $m \equiv 2,3 \pmod 4$,
   then the subgroup $\tScal_{n-1}$ is also a block stabilizer, as 
   $\spinchar{\lambda} = \spinchar{\mu}\upa{\tScal_n}$ for
   the partition $\mu=(m,m-1,\ldots,1)\in\Dcal_{n-1}^-$.
 \end{itemize}
 And in particular,
 \begin{itemize}
  \item[(iii)]
   if $n=6$ then $\spinchar{4,2}$ 
   is the induced of an irreducible character of degree 2
   of $\tScal_{n:2}\cong 3^2\colon Q_8\colon 2$, and if $n=9$ then
   $\spinchar{6,2,1}$ is an induced linear character of
   order six of the subgroup $2\times L_2(8)\colon 3$.
 \end{itemize}
\end{corollary}

In order to prove Theorem~\ref{thm:Sn}, we have to verify the
minimality of the triples given. As this involves information on the
imprimitive faithful characters of $\tAcal_n$, we make use of our
results in Section~\ref{sec:alt}. Note that these are independent of
Theorem~\ref{thm:Sn}.

\begin{proof}[Proof of Theorem~\ref{thm:Sn}]
 With the help of Corollaries~\ref{cor:mainSn} and \ref{cor:mainAn} we
 recursively trace an imprimitive character to a subgroup from which it
 is induced until we find a minimal triple. Let $\lambda\in\Dcal_n^+$.
 If a constituent of $\spinchar{\lambda}\downa{\tAcal_n}$ were
 an imprimitive character of $\tAcal_n$, then $\lambda$ would be
 an odd partition, or one of the exceptional cases would hold by
 Corollary~\ref{cor:mainAn}. Hence the triples of Part (i) of
 Theorem~\ref{thm:Sn} are minimal. If $\lambda$ is as in Case (ii)
 of Corollary~\ref{cor:mainSn}, then $\tScal_{n-1}$ is another
 block stabilizer, and $\spinchar{\lambda}$ is the induced of a
 not-self-associate spin character of this group. The latter cannot
 be imprimitive by Corollary~\ref{cor:mainSn}. Note that if $n=4$,
 then $\mu=(2,1)$ and $\spinchar{\mu}$ is linear. This proves the
 minimality of the triples in Case (ii) of Theorem~\ref{thm:Sn}.
 Lastly, we have to consider the exceptional cases of Part (iii) of
 Corollary~\ref{cor:mainSn}: With \textsf{GAP} it is elementary to
 check that the character of degree 2 of $3^2\colon Q_8\colon 2$ which
 induces to $\spinchar{4,2}$ is in fact an induced linear character of
 $3^2\colon 8$, giving the claimed minimal triple. Likewise, we confirm
 the minimality of the remaining triple.
\end{proof}
It is now immediate which irreducible spin characters of $\tScal_n$ are
monomial.
\begin{corollary}
 Let $\chi\in\Irr_-(\tScal_n)$ be an imprimitive monomial spin
 character.  Then we have
 \begin{itemize}
  \item[(i)] $n=4$ and $\chi=\spinchar{3,1}$, or
  \item[(ii)] $n=6$ and $\chi=\spinchar{4,2}$, or
  \item[(iii)] $n=9$ and $\chi=\spinchar{6,2,1}$.
 \end{itemize}
\end{corollary}
\begin{proof}
 Going through our list of minimal triples $(H,\varphi,\chi)$ in
 Theorem~\ref{thm:Sn}, we determine which $\varphi$ are linear. For
 $n\geq 4$ there are no faithful non-trivial linear characters in
 $\Irr_-(\tAcal_n)$. In Case (ii) of Theorem~\ref{thm:Sn} the partition
 $\mu=(m,m-1,\dots,1)$ of $n-1$ is an odd staircase. Hence, if $m\ge
 3$, then the degree of $\spinchar\mu$ is at least two. Only if $m=2$
 the resulting spin character $\spinchar\mu=\spinchar{2,1}$ is linear.
 The exceptional imprimitive characters of Parts (iii) and (iv) are
 evidently monomial.
\end{proof}

\section{The Alternating Groups}\label{sec:alt}
Using the results of the previous section, we may now determine the
irreducible spin characters of $\tAcal_n$ which are induced spin
characters of subgroups, with the help of Clifford Theory.

In analogy to Section~\ref{sec:symm}, we have two types of maximal
subgroups to consider: On the one hand there are the intersections with
the preimages of the maximal parabolic subgroups of $\tScal_n$, i.e.,
the subgroups $\tAcal_{l,n-l} := \tScal_{l,n-l}\cap\tAcal_n$. On the
other hand there are the subgroups $\tAcal_{n:2} := \tScal_{n:2}\cap
\tAcal_n$.

\subsection{Characters Induced from $\tAcal_{l,n-l}$}
We begin by considering the spin characters of the subgroups
$\tAcal_{l,n-l}$. With the help of Lemma \ref{la:parabolSn} it is now
straightforward to determine which characters induce irreducibly.

\begin{lemma}\label{la:parabolAn}
 A faithful irreducible character $\chi$ of $\tAcal_n$ is
 induced from a subgroup $\tAcal_{l,n-l}$ if and only if
 $l=n-1$, in other words $\tAcal_{l,n-l}=\tAcal_{n-1}$,
 and $\chi = \spinchar\lambda\downa{\tAcal_n}$ with
 $\lambda=(m+1,m-1,m-2,\dots,1)\in \Dcal_n$ for some
 $m\in\NN$ with $m\equiv 0$ or $1\pmod 4$. In this case 
 $\chi$ is the induced of either constituent
 of $\spinchar{\mu}\downa{\tAcal_{n-1}}$ with
 $\mu=(m,m-1,\dots,1)\in\Dcal_{n-1}^+$.
\end{lemma}
\begin{proof} 
 Let $\alpha$ be a spin character of $ \tAcal_{n-l,l}$
 such that $\alpha\upa{\tAcal_n}$ is irreducible. If
 $\alpha\upa{\tScal_n}=\spinchar\lambda$ is irreducible, so is
 $\alpha\upa{\tScal_{l,n-l}}$. Therefore by Lemma~\ref{la:parabolSn}
 we have $l=n-1$ and $\alpha\upa{\tScal_{n-1}} = \spinchar{\mu}$
 for some $\mu \in \Dcal_{n-1}^-$. Hence $\alpha\upa{\tScal_{n-1}}$
 is not self-associate, which is a contradiction. Therefore
 $\alpha\upa{\tScal_n}=\spinchar\lambda+\spinchar\lambda^a$ for a not
 self-associate character $\spinchar\lambda \in \Irr_-(\tScal_n)$.

 Let $l, n-l > 1$. Suppose $\alpha\upa{\tScal_{n-l,l}}$ is
 irreducible. Then Lemma~\ref{la:parabolSn} forces $l=3$ and
 $\alpha\upa{\tScal_{n-3,3}}$ is not self-associate, giving
 again a contradiction. Therefore $\alpha\upa{\tScal_{n-l,l}}
 = \varphi + \varphi^a$ for a not self-associate $\varphi \in
 \Irr_-(\tScal_{n-l,l})$, and $\varphi$ induces to $\spinchar{\lambda}$
 or $\spinchar{\lambda}^a$. By Lemma~\ref{la:parabolSn} this is
 impossible if $l, n-l > 1$.
 
 So let $l=n-1$ and suppose that $\alpha\upa{\tScal_{n-1}}
 = \spinchar{\mu}$ for some $\mu \in \Dcal_{n-1}^+$.
 Hence $\spinchar{\mu}\upa{\tScal_n} = \spinchar{\lambda}
 + \spinchar{\lambda}^a$, so by the Branching Rule we
 obtain $\mu = (m,m-1,\ldots,1)$ where $m\equiv 2,3 \pmod
 4$. On the other hand, if $\alpha\upa{\tScal_{n-1}}$ is
 reducible, i.e., $\alpha\upa{\tScal_{n-1}} = \spinchar{\mu} +
 \spinchar{\mu}^a$ for some $\mu \in \Dcal_n^-$, we obtain again that
 $\spinchar{\mu}\upa{\tScal_n}$ induces to either $\spinchar{\lambda}$
 or $\spinchar{\lambda}$, and thus by Lemma~\ref{la:parabolSn} the
 partition $\lambda$ is even, giving another contradiction.
\end{proof}

\subsection{Characters Induced from $\tAcal_{n:2}$}
Let $n=2m$, where $m\geq 2$. By Clifford Theory, an irreducible character
$\alpha\in\Irr_-(\tAcal_{m,m})$ is either invariant under conjugation
in $\tAcal_{n:2}$ or its inertia subgroup is $\tAcal_{m,m}$. In the
first case it is $\alpha\upa{\tAcal_{n:2}}=\alpha^++\alpha^-$ for
two extensions of $\alpha$ to $\tAcal_{n:2}$. In the second case
$\alpha^0:=\alpha\upa{\tAcal_{n:2}}$ is irreducible.

The following result is an immediate consequence of
Corollary~\ref{la:parabolAn}.

\begin{lemma} \label{la:notinvariantAn}
 For $\alpha^0 \in\Irr_-({\tAcal_{n:2}})$ the induced character
 $\alpha^0\upa{\tAcal_n}$ is always reducible.
\end{lemma}

We now consider the case when the inertia subgroup of
$\alpha\in\Irr_-(\tAcal_{m,m})$ is $\tAcal_{n:2}$, and obtain an
analogue to Corollary~\ref{cor:indpairsSn}.

\begin{lemma} \label{la:invariantAn}
 Let $\alpha\in\Irr_-({\tAcal_{m,m}})$ be invariant under the action of
 $\tAcal_{n:2}$, i.e., $\alpha\upa{\tAcal_{n:2}}=\alpha^++\alpha^-$ for
 two extensions $\alpha^+$, $\alpha^-$ of $\alpha$ to $\tAcal_{n:2}$.
 Then $\alpha\upa{\tScal_{m,m}} \in \Irr_-(\tScal_{m,m})$ is
 invariant in $\tScal_{n:2}$, and we have $\alpha^+\upa{\tScal_n} =
 \alpha^-\upa{\tScal_n}$.
\end{lemma}
\begin{proof}
 First, we note that $\alpha\upa{\tScal_{n:2}}=
 \alpha^+\upa{\tScal_{n:2}}+ \alpha^-\upa{\tScal_{n:2}}$ and both
 summands are self-associate. Assume that $\alpha$ is invariant in
 $\tScal_{m,m}$. Then we have $\alpha\upa{\tScal_{m,m}}=\psi+\psi^a$ for
 some irreducible not self-associate $\psi\in\Irr_-(\tScal_{m,m})$ by
 \cite[Theorem 5.10]{HoffmanHumphreys}. Say $\psi=\spinchar\mu \otimes_z
 \spinchar\nu$ for some $\mu,\nu\in\Dcal_m$ with $\sigma(\mu)\neq
 \sigma(\nu)$. In particular, neither $\psi$ nor $\psi^a$ are
 invariant in $\tScal_{n:2}$. Therefore, both $\psi\upa{\tScal_{n:2}}$
 and $(\psi^a)\upa{\tScal_{n:2}}=(\psi\upa{\tScal_{n:2}})^a$
 are irreducible. By the initial decomposition of
 $\alpha\upa{\tScal_{n:2}}$ both constituents are self-associate,
 so $\psi\upa{\tScal_{n:2}}= (\psi\upa{\tScal_{n:2}})^a$. But by
 Corollary~\ref{cor:actionttauonchars} $\spinchar\mu \otimes_z
 \spinchar\nu$ and either $\spinchar\nu \otimes_z \spinchar\mu$ or
 $(\spinchar\nu \otimes_z \spinchar\mu)^a$ are the only constituents
 of $\psi\upa{\tScal_{n:2}}\downa{\tScal_{m,m}}$, a contradiction.
 Therefore $\alpha$ is not invariant in $\tScal_{m,m}$ and
 $\alpha\upa{\tScal_{m,m}}$ is irreducible and invariant under the
 conjugation action of $\tScal_{n:2}$. As $\alpha^+\upa{\tScal_{n:2}}$
 and $\alpha^-\upa{\tScal_{n:2}}$ are self-associate the equality
 $\alpha^+\upa{\tScal_n}=\alpha^-\upa{\tScal_n}$ follows from Lemma
 \ref{la:indvalues}.
\end{proof}
\begin{corollary}\label{cor:inducedpairsAn}
 Under the hypothesis of Lemma~\ref{la:invariantAn} we have
 that $\alpha^+\upa{\tAcal_n}$ is irreducible if and only if
 $\alpha^-\upa{\tAcal_n}$ is irreducible.
\end{corollary}
\begin{proof}
 Without loss of generality let $\alpha^+\upa{\tAcal_n}$ be irreducible.
 Induction to $\tScal_n$ gives that either $\alpha^+\upa{\tScal_n}$ is
 irreducible, which forces $\alpha^-\upa{\tAcal_n}$ to be irreducible
 too by Lemma \ref{la:invariantAn}, or $\alpha^+\upa{\tScal_n} =
 \spinchar{\lambda} + \spinchar{\lambda}^a$ for some $\lambda \in
 \Dcal_n^-$. Were $\alpha^-\upa{\tAcal_n}$ not irreducible in the
 latter case, it would have two constituents inducing irreducibly
 to either $\spinchar{\lambda}$ or $\spinchar{\lambda}^a$. But then
 $\spinchar{\lambda}$ would be self-associate, giving a contradiction.
\end{proof}
\begin{corollary} \label{cor:imprimitiveAn}
 Let the hypothesis of Lemma~\ref{la:invariantAn} hold. The characters
 $\alpha^+\upa{\tAcal_n}$ and $\alpha^-\upa{\tAcal_n}$ are irreducible
 if and only if $n=6$ and $\alpha\upa{\tScal_{3,3}}=\spinchar{2,1}
 \otimes_z \spinchar{2,1}$. In this case, the characters
 $\alpha^+\upa{\tAcal_6}$ and $\alpha^-\upa{\tAcal_6}$ are the two
 conjugate constituents of $\spinchar{4,2}\downa{\tAcal_6}$.
\end{corollary}
\begin{proof}
 Let $\alpha^+\upa{\tAcal_n}$, and therefore by
 Corollary~\ref{cor:inducedpairsAn} also $\alpha^-\upa{\tAcal_n}$, be
 irreducible.
 By Lemma \ref{la:invariantAn} 
 we have $\alpha^+\upa{\tScal_n}=\alpha^-\upa{\tScal_n}$ and $\alpha$ is a
 constituent of 
 $(\spinchar{\mu}\otimes_z \spinchar{\mu})\downa{\tAcal_{m,m}}$
 for some $\mu \in \Dcal_m$. As
 $\alpha^+\upa{\tScal_n}$ is either irreducible or the sum of an associate
 pair of characters, Lemma~\ref{la:kleinevielfachheit} gives that
 $\mu$ is either $(2,1)$ or $(3,2,1)$.
 But if $\mu=(3,2,1)$, both constituents of
 $(\spinchar{\mu}\otimes_z \spinchar{\mu})\downa{\tAcal_{6,6}}$ 
 are not invariant in $\tAcal_{12:2}$, hence
 induce irreducibly to $\alpha^0$.  This leaves the case $n=6$ and
 $\mu=(2,1)$ in which both constituents of
 $(\spinchar\mu\otimes_z \spinchar{\mu})\downa{\tAcal_{3,3}}$ are invariant in
 $\tAcal_{6:2}$.  In this case the extensions $\alpha^+$ and $\alpha^-$ induce
 to the two conjugate constituents of $\spinchar{4,2} \downa{\tAcal_n}$.
\end{proof} 

In analogy to Corollary~\ref{cor:mainSn}, we summarize the 
conclusions of
Lemma~\ref{la:kleidmanwales}, Lemma~\ref{la:parabolAn}, and
Corollary~\ref{cor:imprimitiveAn} in the following corollary.

\begin{corollary}\label{cor:mainAn}
 Let $\chi$ be an imprimitive spin character of $\tAcal_n$.
 Then
 \begin{itemize}
  \item[(i)] $\chi$ is the restriction to $\tAcal_n$ of the not-self-associate
   spin character 
   $\spinchar{\lambda}$ of $\tScal_n$ with
   $\lambda=(m+1,m-1,m-2,\ldots,1)$ for some $m\equiv 0,1 \pmod 4$. It is
   the induced of either constituent of
   $\spinchar{m,m-1,\ldots,1}\downa{\tAcal_{n-1}}$ and hence
   $\tAcal_{n-1}$ is a block stabilizer.
 \end{itemize}
 Additionally, we have the following exceptions:
 \begin{itemize}
  \item[(ii)] 
   If $n=6$ then $\chi$ is either constituent of
   $\spinchar{4,2}\downa{\tAcal_6}$. Both are induced extensions to
   $\tAcal_{6:2}\cong 3^2\colon 8$ of linear characters of 
   the subgroup $\tAcal_{3,3}\cong 3^2\colon4$ of order four.
  \item[(iii)] If $n=9$ then $\chi$ is either constituent of
   $\spinchar{6,2,1}\downa{\tAcal_9}$. Both are induced linear 
   characters of order six of the subgroup $2\times L_2(8)\colon 3$
 \end{itemize}
\end{corollary}

With the help of Corollary~\ref{cor:mainAn} it is now straightforward to
establish the veracity of the claims made in Theorem \ref{thm:An}.

\begin{proof}[Proof of Theorem 1.2.]
 If $\lambda$ is as in Case~(i) of Corollary~\ref{cor:mainAn}, then
 the staircase partition $(m,m-1,\ldots,1)$ of $n-1$ is an even partition.
 Hence the constituents of $\spinchar{m,m-1,\ldots,1}$ are not
 imprimitive, giving the minimality of the triples in Part~(i) of
 Theorem~\ref{thm:An}. The minimality of the exceptional triples is
 immediate, as the corresponding imprimitive characters are monomial.
\end{proof}

Going through the list of minimal triples $(H,\varphi,\chi)$ in \ref{thm:An}
we can easily identify the monomial characters of $\tAcal_n$.
\begin{corollary}
  Let $\chi\in\Irr_-(\tAcal_n)$ be an imprimitive,
  monomial faithful character. Then
  \begin{itemize}
    \item[(i)] $n=6$ and $\chi$ is a constituent of $\spinchar{4,2}
      \downa{\tAcal_6}$, or
    \item[(ii)] $n=9$ and $\chi$ is a constituent of
      $\spinchar{6,2,1}\downa{\tAcal_{9}}$.
  \end{itemize}
\end{corollary}
\begin{proof}
  The constituents of $\spinchar\mu\downa{\tAcal_{n-1}}$ for
  $\mu=(m,m-1,\dots,1)\in\Dcal_{n-1}^+$ with $m\geq 4$ have degree at least 4.
  On the other hand, we have already mentioned, that the characters $\varphi$
  of the exceptions (ii) and (iii) of Theorem~\ref{thm:An} are linear.
\end{proof}

\subsection{The Sixfold Covers of $\Acal_6$ and $\Acal_7$}
As we have pointed out in the introduction, 
the Schur covers of the alternating groups 
$\Acal_6$ and $\Acal_7$ are sixfold covers (see
\cite[Theorem~2.11]{HoffmanHumphreys})
in contrast to the double covers we have considered so far. Hence to
complete our classification of the imprimitive faithful characters, we
still have to consider the non-split extensions $6\dot{}\Acal_6$ and
$6\dot{}\Acal_7$. For the sake of completeness we will also consider
the triple covers $3\dot{}\Acal_6$ and $3\dot{}\Acal_7$.

The characters and the maximal subgroups of the groups considered are
given in the ATLAS \cite{ATLAS} or in the \textsf{GAP} character table 
library. 
In the sequel, if the character table of a group is given in the ATLAS,
we denote its characters as they are denoted there. Note that, as in
each of the ATLAS tables considered, of a pair of complex conjugate faithful
characters only one member is printed, by a slight abuse of notation, more than 
one character may have the same label. Here we distinguish these characters
by writing $\bar\chi$ for the character whose proxy is $\chi$.

In the case of $6\dot{}\Acal_6$
we deduce from the character degrees that any faithful imprimitive
character is necessarily induced from a maximal subgroup isomorphic to
$3\times \tAcal_5$. Indeed, an analysis of the character
table of this group with $\textsf{GAP}$ shows that 
each faithful character of degree two induces to 
either $\chi_{21}$ or $\chi_{22}$ (or its complex conjugates).
Since there is no maximal subgroup of $3\times
\tAcal_5$ of index two, this yields minimal triples in the sense of
Theorems~\ref{thm:Sn} and \ref{thm:An}.

By the same arguments in the case of $3\dot{}\Acal_6$,
we have two isomorphism types of maximal subgroups to consider: 
For $3\times \Acal_5$ the two
non-trivial linear characters induce irreducibly to give both faithful
characters $\chi_{16}$ and $\bar\chi_{16}$.
Further imprimitive characters
arise by inducing both linear characters of $3\times \Scal_4$ which have
order six. This yields $\chi_{18}$ and $\bar\chi_{18}$
of $3\dot{}\Acal_6$.

Our treatment of $6\dot{}\Acal_7$ and $3\dot{}\Acal_7$ is the same: Again by
examining the character degrees of $6\dot{}\Acal_7$, we conclude that no
faithful character is imprimitive. The situation is somewhat different for
$3\dot{}\Acal_7$, however. Here we have three isomorphism types of maximal
subgroups to consider which contribute to imprimitive characters of
$3\dot{}\Acal_7$: Inducing any faithful character of degree 3 of
$3\dot{}\Acal_6$ gives an irreducible character of degree 21. 
More precisely,
both characters $\chi_{14}$ and $\chi_{15}$ of $3\dot{}\Acal_6$ 
induce to the character $\chi_{21}$ of $3\dot{}\Acal_7$.
This again yields minimal triples, since there is no maximal subgroup of
$3\dot{}\Acal_6$ of index 3.
Also,
the characters $\chi_{20}$ and $\bar\chi_{20}$ of $3\dot{}\Acal_7$ 
are induced linear characters of order six of
$3\times \Scal_5$. And lastly,
the faithful characters $\chi_{19}$ and $\bar\chi_{19}$ of $3\dot{}\Acal_7$ 
are induced non-trivial linear characters of a maximal
subgroup isomorphic to $3\times \mathrm{L}_2(7)$.

We summarize the above in the following theorem. 
\begin{theorem}\label{thm:6An}
 For $3\dot{}\Acal_6$ the minimal triples $(H,\varphi,\chi)$ are as follows:
 \begin{itemize}
  \item[(i)] $H=3\times \tAcal_5$, $\varphi$ is a non-trivial linear
   character, and $\chi$ is either $\chi_{16}$ or $\bar\chi_{16}$.
  \item[(ii)] $H=3\times \Scal_4$, $\varphi$ is a
   linear character of order six, and $\chi$ is $\chi_{18}$
   or $\bar\chi_{18}$.
 \end{itemize}
 For $3\dot{}\Acal_7$ we have
 \begin{itemize}
  \item[(i)] $H=3\dot{}\Acal_6$, $\varphi \in \{\chi_{14},\chi_{15}\}$, 
   and $\chi=\chi_{21}$, or $\varphi \in \{\bar\chi_{14},\bar\chi_{15}\}$
   and $\chi = \bar\chi_{21}$.
  \item[(ii)] $H=3\times \Scal_5$, $\varphi$ is a 
   linear character of order six, and $\chi$ is either $\chi_{20}$ or
   $\bar\chi_{20}$.
  \item[(iii)] $H=3\times L_2(7)$, $\varphi$ is a
   non-trivial linear character, and $\chi$ is $\chi_{19}$ or
   $\bar\chi_{19}$.
 \end{itemize}
 For $6\dot{}\Acal_6$ we have
 \begin{itemize}
  \item[(i)] $H=3\times \tAcal_5$, $\varphi$ has degree two, and $\chi$ is
   either $\chi_{21}$ or $\chi_{22}$ or one of their complex conjugates.
 \end{itemize}
 None of the faithful ordinary characters of $6\dot{}\Acal_7$ are imprimitive.
 \label{thm:odd}
\end{theorem}
\subsection*{Acknowledgment} The authors are indebted to Kay Magaard and
Gerhard Hiss for many helpful discussions on the subject.
\bibliographystyle{jcm}

\end{document}